\documentclass[a4paper,12pt]{amsart}
\usepackage[labelsep=colon]{caption}
\usepackage{cite}
\usepackage{float}
\restylefloat{figure}
\usepackage{graphicx}
\usepackage[colorlinks=true]{hyperref}
\hypersetup{allcolors=[rgb]{0.9,0.1,0.4}}
\usepackage{mathrsfs}
\usepackage[hang]{subfigure}
\usepackage{tikz}
\usetikzlibrary{arrows,%
  decorations.markings,%
  decorations.pathmorphing,%
  matrix,%
  shapes}
\usepackage{wrapfig}

\hypersetup{
  pdftitle={A Mathematical Approach to the Hierarchical Structure of
    Languages},%
  pdfauthor={Nils A. Baas},%
}

\tikzstyle{midarrow} = [decoration={markings,mark=at position .5 with
  \arrow{>}}, postaction={decorate}]

\newcommand{\boks}{\mathop{}\!\square\mathop{}\!}
\newcommand{\diagram}[3]{\matrix (#1) [matrix of math nodes,row
  sep={#2},column sep={#3},text height=1.5ex,text depth=0.25ex]}
\renewcommand{\H}{\mathscr{H}}
\renewcommand{\L}{\mathscr{L}}
\renewcommand{\P}{\mathscr{P}}

\newcommand{\sets}{\mathrm{Sets}}

\title{A Mathematical Approach to the Hierarchical Structure of
    Languages}
\author[N.A. Baas]{Nils A.\ Baas}
\date{May 2, 2018}
\address{Department of Mathematical Sciences, NTNU, NO-7491 Trondheim,
  Norway}
\email{nils.baas@ntnu.no}

\begin{document}
\maketitle


\section{Introduction}

We will in this paper suggest how the mathematical concept of
hyperstructures may be a useful tool in the study of the higher,
hierachical structure of languages.

By a language $\L$ we will mean a spoken, written, formal or geometric
language.  We will assume that in the context of a language there are
some basic units: sounds, letters, geometric forms, signs, etc.  They
represent the ``atoms'' of the language.  We suggest that
Hyperstructures \cite{Cog, HAM, NSCS, NS, SO, HOA, HOS, MHS, PHS}
represent a useful framework in which to represent a language
structure, and suggest the following procedure.

\section{Hyperstructures}

Let us denote the collection or set of basic linguistic units by
$X_0$.  Let $\P(X_0)$ stand for all (finite) subsets or a
subcollection.  We will then have an assignment (technically a
presheaf in many cases):
\begin{equation*}
  \Omega_0 \colon \P(X_0) \to \sets,
\end{equation*}
for example ``giving some meaning'' to a collection of sounds (or
letters).  This is what we call the Observation part (property,
state).

Next we consider collections of basic units with properties
\begin{equation*}
  (S_0,\omega_0)
\end{equation*}
and we want to bind these together to meaningful entities
\begin{equation*}
  \Gamma_1 = \{(S_0,\omega_0) \mid S_0 \subseteq X_0, \omega_0 \in
  \Omega_0(S_0)\}.
\end{equation*}

$B_0 \colon \Gamma_1 \to \sets$ is another assignment (presheaf) such
that $B_0(S_0,\omega_0)$ is the set of bonds, ways of gluing the
element in $S_0$ together to meaningful units, like words.

This is the process we will iterate in order to get a
hyperstructure. We think of
\begin{itemize}
\item[] $\Omega_0$ as the semantic assignment, ``meaning''
\end{itemize}
and
\begin{itemize}
\item[] $B_0$ as the syntactic (``grammatical'') structure assignment.
  Admissible generation of words.
\end{itemize}

In hyperstructures we can also define composition of bonds, by a
``gluing'' process (see \cite{HOA}) which is a huge
generalization of concatenation of letters for example.  From bonds
$b_1,b_2,\ldots,b_k$ we can generate new bonds:
\begin{equation*}
  b = b_1 \boks b_2 \boks \cdots \boks b_k.
\end{equation*}
These assignments give a ``grammar'' both for syntax and semantics.
For example Chomsky grammars and their generative rules will be
covered by bond composition in the form of concatenation. We may also
consider higher dimensional geometric alphabets and higher dimensional
gluing. The semantic part may also play a role in the composition.

Now we can proceed to produce the higher tiers of the language:

Put
\begin{equation*}
  X_1 = \{ b_0 \mid b_0 \in B_0(S_0,\omega_0)\}.
\end{equation*}
Then we choose or construct (a ``presheaf''), specific for each
language:
\begin{equation*}
  \Omega_1 \colon \P(X_1) \to \sets
\end{equation*}
from
\begin{equation*}
  \Gamma_1 = \{(S_1,\omega_1)\}
\end{equation*}
and
\begin{equation*}
  B_1 \colon \Gamma_1 \to \sets.
\end{equation*}
This means that $B_1(S_1,\omega_1)$ represent bonds of bonds, for
example words of words bound together in meaningful admissable
sentences given by the bond presheaves.

Then we can construct arbitrarily many tiers getting a formal
hyperstructure $\H(\L)$. 

We also have ``boundary'' maps connecting the tiers:
\begin{equation*}
  \partial_i \colon X_{i + 1} \to \P(X_i)
\end{equation*}
which dissolve the bonds --- or put differently: ``a bond knows what
it binds'' or ``a sentence knows its words'' (linguistic version).

Starting with a top bond and iteratively applying the $\partial_i$'s
is like given a sentence and constructing a parsing tree:
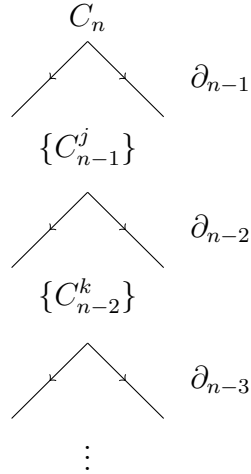
\begin{center}
  \begin{tikzpicture}
    \begin{scope}
      \draw[midarrow] (0,0) -- (1,-1);
      \draw[midarrow] (0,0) -- (-1,-1);
      \node[above] at (0,0){$C_n$};
      \node[below] at (0,-1){$\{C_{n - 1}^j\}$};
      \node at (1.75,-0.5){$\partial_{n - 1}$};
    \end{scope}

    \begin{scope}[yshift=-2cm]
      \draw[midarrow] (0,0) -- (1,-1);
      \draw[midarrow] (0,0) -- (-1,-1);
      \node[below] at (0,-1){$\{C_{n - 2}^k\}$};
      \node at (1.75,-0.5){$\partial_{n - 2}$};
    \end{scope}

    \begin{scope}[yshift=-4cm]
      \draw[midarrow] (0,0) -- (1,-1);
      \draw[midarrow] (0,0) -- (-1,-1);
      \node[below] at (0,-1){$\vdots$};
      \node at (1.75,-0.5){$\partial_{n - 3}$};
    \end{scope}
  \end{tikzpicture}
\end{center}
i.e., describing $n$-levels of generalized ``morphemes''.

\emph{We suggest that this is the basic general structure of language,
and that any reasonable language is governed by a Hyperstructure, and
it may be very useful to identify it.}

Hyperstructures are really organizing tools, and so are languages for
communication, thought and learning.

\section{Consequences}

Notice that even starting with a finite set of basic units, infinite
numbers of bonds and meanings can be generated for any tier.  Our
approach covers semiotic systems as well.

Another advantage with hyperstructures is the possibility of passing
from local to global via an extended Grothendieck type topology (see
\cite{HOA} where this is introduced).

This means that we consider tiers of bonds
\begin{center}
  \begin{tikzpicture}
    \diagram{d}{1em}{2.5em}{
      B_n & \sets\\
      B_{n - 1} & \sets\\
      \vdots &\\
      B_0 & \sets\\
    };

    \path[->,font=\scriptsize,above,midway]
      (d-1-1.east)
      edge[decorate,decoration={snake,amplitude=.4mm,segment
        length=3mm,post length=1mm}] node{$\Lambda_n$} (d-1-2.west)
      (d-2-1.east)
      edge[decorate,decoration={snake,amplitude=.4mm,segment
        length=3mm,post length=1mm}] node{$\Lambda_{n - 1}$}
        (d-2-2.west)
      (d-4-1.east)
      edge[decorate,decoration={snake,amplitude=.4mm,segment
        length=3mm,post length=1mm}] node{$\Lambda_0$} (d-4-2.west);
  \end{tikzpicture}
\end{center}
where the $\Lambda_i$'s are ``meaning'' functions or assignments at
each level. We have developed rules for when we from these
\emph{levelwise} meanings can deduce a global meaning.  This is done
by what we call a globalizer (generalizing sections in sheaves), and
shows how global meanings can be integrated and also that in many
situations you cannot just pass from local to global.  Tiers or levels
are needed.  This extends the linguistic notion of \emph{``duality of
  patterning''}.  But duality is just for two levels, in general many
levels are needed in an efficient language.  Our approach is
\emph{``multiplicity of patterning''}!

In general: Data, Stimuli, Inputs, Outputs may be unstructured ---
have no immediate meaning.

Hence put a Hyperstructure on the situation (which may be inherent in
the living brain or a synthetic one) and via a suitable ``globalizer''
structure and meaning will appear (unless random!).

This is a way to handle the ``Generalized Binding Problem'' (see
\cite{HOS}).  If two languages are $\H$-structured communication
will (should) take place in an $\H$-structured way and understanding
represents some kind of structural compatibility or resonance.

Finally, these structures also relate to memory and learning.  For
example a \emph{machine} organized as an $\H$-structure will take data
(stimuli), process them and adjust according to a global goal state or
meaning.

This is \emph{automatic} learning.  Learning in general requires or at
least profits from a Hyperstructure.  Hyperstructured libraries of
events or data are important in learning and memory and
$\H$-structures facilitates the process.

    


We hope that this note will initiate a further study of
hyperstructures in languages.

\subsection*{Acknowledgements}
I would like to thank M.~Thaule for his kind technical assistance in
preparing the manuscript.

\end{document}